\documentclass{article}
\usepackage{amsmath}
\usepackage{graphicx}
\usepackage[all]{xypic}
\usepackage{latexcad}
\usepackage{amsfonts}
\usepackage{amssymb}
\makeatletter
\newtheorem{thm}{Theorem}
\newtheorem{prop}[thm]{Proposition}
\newtheorem{defn}[thm]{Definition}
\newtheorem{rem}[thm]{Remark}
\newtheorem{ex}[thm]{Example}

\makeatother

\begin{document}

\title{Cospans and spans of graphs: a categorical algebra for the sequential and
parallel composition of discrete systems}
\author{L. de Francesco Albasini, N. Sabadini, R.F.C. Walters\\Dipartimento di Informatica e Comunicazione,\\Universit\`{a} dell' Insubria, Varese, Italy}
\maketitle
\date{}
\begin{abstract}
We develop further the algebra of cospans and spans of graphs introduced by
Katis, Sabadini and Walters \cite{KSW00biwim} for the sequential and parallel
composition of processes, adding here data types.
\end{abstract}

\section{Introduction}

This paper develops further the algebra for the sequential and parallel
composition of systems introduced in the two papers \cite{KSW97bspan graph},
\cite{KSW00biwim}. Whereas those papers dealt with the finite state control,
here we add data structures. As in \cite{KSW00biwim} the sequential
composition is a cospan composition, the parallel a span composition.

The plan of the paper is as follows. We begin with the most abstract notion of
a system with sequential and parallel interfaces. In section 3 we
make simplifying assumptions arriving in the section 4 with an algebra
that is, in effect, an implementable programming language of systems. The
reader should be aware that the word \emph{system} has an increasingly
specific meaning in successive sections of the paper. The motive for
proceeding in this way is to show that what may appear as arbitrary and
unmotivated in section 4 actually arises in a natural way from general
considerations. In addition, for some applications a different set of
simplifying assumptions may be more appropriate. As examples of programs 
in the language we indicate in section 5 how sequential programming, 
classical concurrency examples, hierarchy and change of geometry may be expressed.

An important element of this paper is the matrix calculus which arises from
the fact that categories of spans  in an extensive category \cite{CLW93carbonilackwalters}
 have direct sums. 
It allows an explicit
relation between programs with data types, and finite automata which express
the control structure of the program.

Another important element is the role of the distributive law in various
roles, including flattening hierarchy.

The work has been influenced by our study in \cite{dFASW08dmarkov},
\cite{dFASW08equantum} of probabilistic and quantum automata.

In this paper we concentrate on the operations of the algebra, and its
expressivity, rather than the equations satisfied.

\section{Systems with sequential and parallel interfaces}

We represent systems by (possibly infinite) graphs of states and transitions, to which
we will be adding extra struture.
By a graph $G$ we mean here a set of states $states(G)$ a set of transitions
$transitions(G)$ and two functions $source,target:transitions(G)\rightarrow states(G)$
which specify the source state and target state of a transition.

\subsection{Sequential interfaces}

In order to compose sequentially one system with another both systems must
have appropriate interfaces. The idea comes from the sequential composition of
automata, which occurs for example in Kleene's theorem: certain states (final
states) of one automata are identified with certain states (initial states) of
another. Here we replace initial and final states by graph morphisms into the
graph of the system.

\begin{defn}
A \emph{system with sequential interfaces} is a cospan $\gamma_{0}%
:A\rightarrow G\leftarrow B:\gamma_{1\text{ }}$of graphs. The graph $G$ is the
graph of the system; $A$, and $B$ are the graphs of the interfaces. We write
this also as $(G,\gamma_{0},\gamma_{1}):A\rightarrow B$ or even just as
$G_{B}^{A}.$ Composition of systems is by pushout. The category of systems
with sequential interface is $Cospan(Graph).$ A behaviour of $G_{B}^{A}$ is a
path in the central graph $G$.
\end{defn}

Notice that in speaking of the \emph{category} of cospans we should consider
cospans only up to an isomorphism of the central graph of the cospan. In
practice we will always consider representative cospans, and any equation we
state will be true only up to isomorphism. The same proviso should be applied
to our discussion later of spans, and systems.

\subsection{Parallel interfaces}

Similarly, to compose in communicating parallel two systems each system must
have a parallel interface. The idea here comes, for example, from circuits. A
circuit component has a physical boundary and transitions of the circuit
component produce transitions on the physical boundary. Joining two circuit
components, the transitions of the resulting system are restricted by the fact
that the transitions on the common boundary must be equal. We describe the
relation between transitions of the system $G$ and the transitions on a
boundary $X$ by a graph morphism $G\rightarrow X$. To obtain a category when
we compose we require that a system has two parallel interfaces.

\begin{defn}
A \emph{system with parallel interfaces} is a span $\partial_{0}:X\leftarrow
G\rightarrow Y:\partial_{1\text{ }}$of graphs. The graph $G$ is the graph of
the system; $X$, and $Y$ are the graphs of the interfaces. We write this also
as $(G,\partial_{0},\partial_{1}):X\rightarrow Y$ or even just as $G_{X,Y}.$
Composition of systems is by pullback. The category of systems with sequential
interface is $Span(Graph).$ A behaviour of $G_{X,Y}$ is a path in the central
graph $G$.
\end{defn}

\subsection{ Combined sequential and parallel interfaces}

\begin{defn}
A \emph{system with sequential and parallel interfaces} consists of a
commutative diagram of graphs and graph morphisms%
\[
\xymatrix{
G_0\ar[d]_{\gamma_0}  & X \ar[l]_{\partial_0}\ar[d]_{\gamma_0}\ar[r]^{\partial
_1}&G_1\ar[d]_{\gamma_0}\\
A & G\ar[l]_{\partial_0} \ar[r]^{\partial_1}&B\\
G_2\ar[u] ^{\gamma_1}& Y\ar[l]_{\partial_0}\ar[u]^{\gamma_1}\ar[r]^{\partial
_1}&G_3\ar[u]^{\gamma_1}}%
\]
or more briefly, when we are not emphasizing the corner graphs, as%
\[
\xymatrix{
\bullet\ar[d] & X \ar[l]\ar[d]\ar[r]&\bullet\ar[d]\\
A & G\ar[l] \ar[r]&B\\
\bullet\ar[u] & Y\ar[l]\ar[u]\ar[r]&\bullet\ar[u]
}%
\]
We denote such a system very briefly as $G_{Y;A,B}^{X}$, or even $G_{Y}^{X}$
or $G_{A,B}$ or even just $G$, depending on the context. A behaviour of
$G_{Y;A,B}^{X}$ is a path in the central graph $G$. Another useful notation is
as follows: given an object $O$ in the diagram we denote the four adjacent
objects by $O^{\rightarrow},O^{\leftarrow},O^{\downarrow}$ and $O^{\uparrow};$
for example, $G^{\leftarrow\uparrow}=G^{\uparrow\leftarrow}=G_{0}.$
\end{defn}

Such a system may be regarded in two ways:\ (i) as three systems with parallel
interfaces, the first $X_{G_{0},G_{1}}$ and third $Y_{G_{2},G_{3}}$ being
sequential interfaces to the second $G_{A,B}$; or (ii) as three systems with
sequential interfaces, two ($A_{G_{2}}^{G_{0}},$ $B_{G_{3}}^{G_{1}}$) being
parallel interfaces to the other ($G_{Y}^{X}$). The point is that to compose in 
parallel a system with sequential interfaces requires that the sequential interfaces also have
parallel interfaces. It is not necessary that the parallel interfaces themselves
have parallel interfaces, since interfaces are identified, not composed,
 in the composition. A similar remark applies to sequential
composition.
Notice that for simplicity we
have used the same symbols $\gamma_{0},\gamma_{1}$ for all the sequential
interface morphisms and similarly $\partial_{0},\partial_{1}$ for all the
parallel interface morphisms.

\subsubsection{Operations on systems}

\begin{defn}
Two systems $G_{Y;A,B}^{X}$ and $H_{W;B,C}^{Z}$ admit a compositions by
pullback, the \emph{parallel }(or horizontal) \emph{composition}, denoted
$G_{Y;A,B}^{X}||H_{W;B,C}^{Z}.$
\end{defn}

Of course, certain corner graphs of $G$ and $H$ are required to be the same.
This applies also in the next definition.

\begin{defn}
Two systems $G_{Y;A,B}^{X}$ and $K_{Z;D,E}^{Y}$ admit a compositions by
pushout, the \emph{sequential} (or vertical) \emph{composition}, denoted
$G_{Y;A,B}^{X}\circ K_{Z;D,E}^{Y}.$
\end{defn}

\begin{rem}
Given four systems $G_{Y;A,B}^{X},$ $H_{V;B,C}^{U}$, $K_{Z;D,E}^{Y}$,
$L_{W;E,F}^{V}$ \ in the following configuration
%
%
%
%
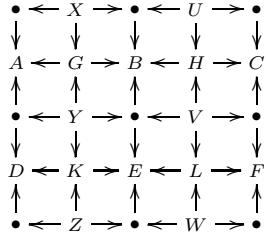
\[
\def\objectstyle{\scriptstyle}
\def\labelstyle{\scriptstyle}
\xymatrix@-1.2pc{\bullet\ar[d]&X\ar[l]\ar[d]\ar[r]&
\bullet\ar[d]&U\ar[l]\ar[d]\ar[r]&\bullet\ar[d]\\
A&G\ar[l]\ar[r]&B&H\ar[l]\ar[r]&C\\
\bullet\ar[u]\ar[d]&Y\ar[l]
\ar[u]\ar[d]\ar[r]&\bullet\ar[u]\ar[d]&V\ar[l]\ar[r]\ar[u]\ar[d]&\bullet\ar[u]\ar[d]\\
D&K\ar[l]\ar[r]&E&L\ar[l]\ar[r]&F\\
\bullet\ar[u]&Z\ar[l]\ar[u]\ar[r]&\bullet\ar[u]&W\ar[l]\ar[u]\ar[r]&\bullet\ar[u]}%
\]
%
%
%
%
%
%
%
%
%
%
%
%
%
there is a comparison map%
\[
(G_{Y;A,B}^{X}||H_{V;B,C}^{U})\circ(K_{Z;D,E}^{Y}||L_{W;E,F}^{V}%
)\rightarrow(G_{Y;A,B}^{X}\circ K_{Z;D,E}^{Y})||(H_{V;B,C}^{U}\circ
L_{W;E,F}^{V}),
\]
satisfying appropriate (lax monoidal) coherence equations, which however is \emph{not in
general} an isomorphism. This reflects the fact that the left-hand expression
involves more synchronization than the right.
\end{rem}

\begin{defn}
The \emph{product} $G_{Y;A,B}^{X}\times H_{W;C,D}^{Z}$ of two systems
$G_{Y;A,B}^{X},$ $H_{W;C,D}^{Z}$ is formed by taking the product of all the
objects and arrows in $G$ with the corresponding objects and arrows in the
$H$; briefly
%
%
%
%
%
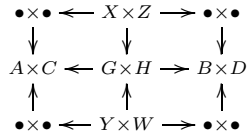
\[
\def\objectstyle{\scriptstyle}
\def\labelstyle{\scriptstyle}
\xymatrix @-1.2pc{
\bullet\times\bullet\ar[d]  & X\times Z \ar[l]\ar[d]\ar[r]&\bullet
\times\bullet\ar[d]\\
A \times C& G\times H\ar[l]\ar[r]&B\times D\\
\bullet\times\bullet\ar[u] & Y\times W\ar[l]\ar[u]\ar[r]&\bullet\times
\bullet\ar[u]}%
\]
%
%
%
%
%
%
%
%
%
%
%
%
%
\end{defn}

\begin{defn}
The \emph{sum} $G_{Y;A,B}^{X}\boxplus H_{W;C,D}^{Z}$ of two systems
$G_{Y;A,B}^{X},$ $H_{W;C,D}^{Z}$ is formed by taking the sum of all the
objects and arrows in $G$ with the corresponding objects and arrows in the
$H$; briefly
%
%
%
%
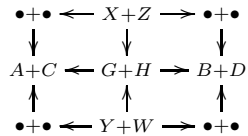
\[
\def\objectstyle{\scriptstyle}
\def\labelstyle{\scriptstyle}
\xymatrix @-1.2pc{
\bullet+\bullet\ar[d]  & X+ Z \ar[l]\ar[d]\ar[r]&\bullet
+\bullet\ar[d]\\
A + C& G+ H\ar[l]\ar[r]&B+ D\\
\bullet+\bullet\ar[u] & Y+ W\ar[l]\ar[u]\ar[r]&\bullet+
\bullet\ar[u]}%
\]
%
%
%
%
%
%
%
%
%
%
%
%
%
\end{defn}

The last part of the algebra of systems consists of a number of constants.

\begin{defn}
The constants of the algebra are systems constructed from the constants of the
distributive category structure of $Sets$ \cite{W89mybook},\cite{L66Lawevere
sym logic},\cite{E75elgot}.
\end{defn}

When we describe in a later section a programming language there will be of
course also as constants the operations of data types; the particular language
we describe has the natural numbers together with predecessor and successor.

\section{Simplifying Assumptions}

We introduce a number of simplifying assumptions with the aim of arriving at a
implementable programming language for systems. As we do this we will be
considering also certain important derived operations of the algebra.

\subsection{Simplifying the interfaces}

\bigskip

\noindent\textbf{Assumption 1. }\textit{We assume from now on that in a system
with sequential and parallel interfaces }$G$\textit{ as described above the
corner graphs }$G^{\leftarrow\uparrow},G^{\rightarrow\uparrow},G^{\leftarrow
\downarrow},G^{\rightarrow\downarrow}$\textit{ each have \emph{one state and
no transitions}, that the graphs }$A,B$\textit{ each have \emph{one state},
and that the graphs }$X,Y$\textit{ have \emph{no transitions}.} \bigskip

The idea is that in many cases the sequential interface consists only of
states with no transitions, whereas the parallel interfaces are ``stateless'',
that is, consist of transitions and one state. The assumptions are appropriate
for message passing communication but not for systems in which there is
communication by shared variables, since this requires that the parallel
interfaces have state. It is not difficult to make assumptions for this type
of communication but we prefer here to make the simpler assumption.

Given the assumption we may ignore the corner graphs of a system so that it
consists of five graphs $G,A,B,X,Y$ and the four graph morphisms%
\[
\xymatrix{
& X \ar[d]_{\gamma_0}\\
A & G\ar[l]_{\partial_0} \ar[r]^{\partial_1}&B\\
& Y\ar[u]^{\gamma_1}&}.
\]
Since the single states of $A$ and $B$ need not have a name, we may sometimes
confuse $A$ and $B$ with $transitions(A)$ and $transitions(B)$ respectively.
We may think of $A,B,X,Y$ as sets, and of $A$ and $B$ as labels for the
transitions in $G$ (the graph morphisms $\partial_{0}$, $\partial_{1}$
providing the labelling).

As a consequence of the simple form of the corner graphs of a system we have
the following result.

\begin{prop}
The \emph{parallel composite} $G_{Y;A,B}^{X}||H_{W;B,C}^{Z}$ \ of two systems
has top sequential interface $X\times Z$, and bottom sequential interface
$Y\times W$; we can summarize this by the formula $G_{Y;A,B}^{X}%
||H_{W;B,C}^{Z}=(G||H)_{Y\times W;A,C}^{X\times Z}$. The \emph{sequential
composite} $G_{Y;A,B}^{X}\circ H_{Z;C,D}^{Y}$ \ has left parallel interface
the graph with one vertex and transitions $transitions(A)+transitions(C)$
which we denote with some \emph{abuse of notation} as $A+C$, and similarly
right parallel interface $B+D$; we can summarize this by the formula
$G_{Y;A,B}^{X}\circ H_{Z;C,D}^{Y}=(G\circ H)_{Z;(A+C),(B+D)}^{X}$. Trivially,
$G_{Y;A,B}^{X}\times H_{W;C,D}^{Z}=(G\times H)_{Y\times W;A\times
C,B\times D}^{X\times Z}$.
\end{prop}

Notice that the class of systems we are considering is closed under sequential
and parallel composition and product, but is not closed under the operation of
sum since the resulting system will have parallel interfaces with two states,
not one.

We now introduce two derived operations similar to the sequential composite
and the sum, but which are local in the sense that the parallel interfaces are
fixed. Intuitively they are sequential operations within a fixed parallel protocol.

\begin{defn}
The \emph{local sequential composition} $G_{Y;A,B}^{X}\bullet H_{Z;A,B}^{Y}$
of two systems $G_{Y;A,B}^{X},$ $H_{Z;A,B}^{Y},$ is formed from $G_{Y;A,B}%
^{X}\circ H_{Z;A,B}^{Y}$ by composing with appropriate codiagonals as
follows:
%
%
%
%
\[
\def\objectstyle{\scriptstyle}
\def\labelstyle{\scriptstyle}
\xymatrix @-1.2pc{
\bullet\ar[d]&
\bullet\ar[d]\ar[l]_{1}\ar[r]^{1}&
\bullet\ar[d]&
X\ar[l]\ar[d]\ar[r]&
\bullet\ar[d]&
\bullet\ar[l]_{1}\ar[d]\ar[r]^{1}&
\bullet\ar[d]\\
A&
A+A\ar[l]_{\nabla}\ar[r]^{1}&
A + A&
G\circ H\ar[l]\ar[r]&
B+ B&
B+B\ar[l]_{1}\ar[r]^{\nabla}&
B\\
\bullet\ar[u]&
\bullet\ar[l]_{1}\ar[r]^{1}\ar[u]&
\bullet\ar[u] &
 Z\ar[l]\ar[u]\ar[r]&
\bullet\ar[u]&
\bullet\ar[l]_{1}\ar[r]^{1}\ar[u]&
\bullet\ar[u]}%
\]
%
%
%
%
%
%
%
%
%
%
%
%
where the codiagonals $\nabla:A+A\rightarrow A,\nabla:B+B\rightarrow B$ are
codiagonals on transitions, but the identity on the single state.
\end{defn}

\begin{defn}
The \emph{local sum} $G_{Y;A,B}^{X}+ H_{W;A,B}^{Z}$ of two systems
$G_{Y;A,B}^{X},$ $H_{W;C,D}^{Z},$ is formed from $G_{Y;A,B}^{X}\boxplus
H_{W;C,D}^{Z}$ by composing with appropriate codiagonals as follows:
%
%
%
%
\[
\def\objectstyle{\scriptstyle}
\def\labelstyle{\scriptstyle}
\xymatrix @-1.2pc{
\bullet\ar[d]&
\bullet+\bullet\ar[d]\ar[l]_{\nabla}\ar[r]^{1}&
\bullet+\bullet\ar[d]&
X+ Z \ar[l]\ar[d]\ar[r]&
\bullet+\bullet\ar[d]&
\bullet+\bullet\ar[l]_{1}\ar[d]\ar[r]^{\nabla}&
\bullet\ar[d]\\
A&
A+A\ar[l]_{\nabla}\ar[r]^{1}&
A + A&
G+ H\ar[l]\ar[r]&
B+ B&
B+B\ar[l]_{1}\ar[r]^{\nabla}&
B\\
\bullet\ar[u]&
\bullet+\bullet\ar[l]_{\nabla}\ar[r]^{1}\ar[u]&
\bullet+\bullet\ar[u] &
 Y+ W\ar[l]\ar[u]\ar[r]&
\bullet+\bullet\ar[u]&
\bullet+\bullet\ar[l]_{1}\ar[r]^{\nabla}\ar[u]&
\bullet\ar[u]}%
\]
%
%
%
%
%
%
%
%
%
%
%
%
Clearly, 
$G_{Y;A,B}^{X}+ H_{W;A,B}^{Z}=(G+ H)_{Y+W;A,B}^{X+Z}$.
\end{defn}

Now the class of systems we are now considering is closed under the operations
of parallel and sequential composition, product, local sequential and local sum.

\subsection{Finiteness assumptions}

In general, pushouts and pullbacks of infinite graphs are not implementable.
We need to make some finiteness assumptions.\bigskip

\noindent\textbf{Assumption 2. }\textit{We assume that in system }%
$G_{B;X,Y}^{A}$\textit{ that }$A$\textit{ and }$B$\textit{ have a \emph{finite
number of transitions}.\bigskip\ }

This assumption means that the pullbacks in the parallel composition are
implementable. A further consequence of this assumption is that the
transitions of the graph $G$ decompose as a disjoint union
\[
transitions(G)={\bigvee}_{a\in A,b\in B}transitions(G)_{a,b}%
\]
where $transitions(G)_{a,b}$ is the set of transitions labelled by $a\in
A,b\in B$. Denote by $G_{a,b}$ the graph with the same states as $G$ but with
transitions $transitions(G)_{a,b}$.\bigskip

The next assumption will have the effect that our systems have a finite state
automata as control structure. Usually finite state automata are presented as
recognizers of regular languages \cite{HMU01Hopcroftullman}. However the
original work of McCulloch and Pitts \cite{MP43mccullochpitts} introduced
automata as systems with thresholds, that is systems with infinite state
spaces which decomposed into finite sums. Our finiteness assumptions are of
this nature.\bigskip

\noindent\textbf{Assumption 3. }\textit{We assume that the set of states of
the graphs }$G,$ $G^{\uparrow},$ $G_{\downarrow}$ \textit{are given as a
finite disjoint sums}$: states(G)=U_{1}+U_{2}+\cdots+U_{m},$
$states(G^{\uparrow})=X_{1}+X_{2}+\cdots+X_{k},$ $states(G_{\downarrow}%
)=Y_{1}+Y_{2}+\cdots+Y_{l}.$\bigskip

The first effect of this is that each of the graphs $G_{a,b}$ ($a\in A,b\in
B$) breaks up as a matrix of spans of sets.

To see this notice that a graph $G$ is just an endomorphism in $Span(Sets)$.
Further the category $Span(Sets)$ has direct sums, the direct sum of $U$ and
$V$ being $U+V$ with injections the functions $i_{X}:U\rightarrow U+V$,
$i_{Y}:V\rightarrow U+V$ considered as spans. and projections the same
functions but now considered as the opposite spans $i_{X}^{op}:U+V\rightarrow
U$, $i_{Y}^{op}:U+V\rightarrow V$. The commutative monoid structure on
$Span(Sets)(U,V)$ is given by sum and the empty span. Since a graph is just an
endomorphism in $Span(Sets)$ a graph $G$ whose state set is $U+V$ may be
represented as a $2\times2$ matrix of spans $\left(
\begin{array}
[c]{cc}%
G_{U,U} & G_{U,V}\\
G_{V,U} & G_{V,V}%
\end{array}
\right)  $. where for example $G_{U,V}=i_{V}^{op}Gi_{U}$. Further
$G=i_{U}G_{U,U}i_{U}^{op}+i_{V}G_{U,V}i_{U}^{op}+i_{U}G_{V,U}i_{V}^{op}%
+i_{V}G_{V,V}i_{V}^{op}$.

Generalizing this to the case in which the states break up into a disjoint sum
of $n$ subsets Assumption 3 implies that each of the graphs $G_{a,b}$ may be
represented as a $k\times k$ matrix of spans, the $i,j$th entry of which we
will denote $G_{a,b,U_{i},U_{j}}$, or even $G_{a,b,i,j}$. It has a simple meaning: $G_{a,b,U_{i}%
,U_{j}}$\emph{ is the set of transitions of }$G$\emph{ labelled }$a,b$\emph{
whose sources lie in }$U_{i}$\emph{ and whose targets lie in }$U_{j\text{.}}$
\emph{The projections of the span }$G_{a,b,U_{i},U_{j}}$\emph{ are the
projections onto the sources and targets.}

\medskip
It is easy also to expand the matrix to include the functions $\gamma
_{0}:X\rightarrow G,\gamma_{1}:X\rightarrow G$. The resulting matrix has
columns indexed by $X_{1},X_{2},\cdots,X_{k},U_{1},U_{2},\cdots,U_{l}$ and
rows indexed by $Y_{1},Y_{2},\cdots,Y_{l},U_{1},U_{2},\cdots,U_{m}$; as an
example when $k=l=m=2$ the matrix has the form%
\[%
\begin{array}
[c]{c||cc|cc|}%
G_{a,b} & X_{1} & X_{2} & U_{1} & U_{2}\\\hline\hline
Y_{1} & 0 & 0 & G_{a,b,U_{1},Y_{1}} & G_{a,b,U_{2},Y_{1}}\\
Y_{2} & 0 & 0 & G_{a,b,U_{1},Y_{2}} & G_{a,b,U_{2},Y_{2}}\\\hline
U_{1} & G_{a,b,X_{1},U_{1}} & G_{a,b,X_{2},U_{1}} & G_{a,b,U_{1},U_{1}} &
G_{a,b,U_{2},U_{1}}\\
U_{2} & G_{a,b,X_{1},U_{1}} & G_{a,b,X_{2},U_{2}} & G_{a,b,U_{1},U_{2}} &
G_{a,b,U_{2},U_{2}}\\\hline
\end{array}
\]
where $0$ denotes the empty span.

\begin{ex}
The function $predecessor:N\rightarrow N+1$ which returns an error if the
argument is $0$ but otherwise decrements, may be considered as a system with
trivial parallel interfaces, top sequential interface $N$ bottom sequential
interface $N+1$ and central graph having states $N+N+1$, transitions $N$ and
$source:N\rightarrow N+(N+1)=inj_{N}$, $\ target:N\rightarrow
N+(N+1)=inj_{(N+1)}\cdot predecessor$. (This is the usual picture of a function as a
graph on the disjoint union of the domain and codomain, with edges relating
domain elements and their images,) \ We call this system $pred.$ The matrix is%
\[%
\begin{array}
[c]{c||c|ccc|}%
pred & N & N & N & 1\\\hline\hline
N & 0 & 0 & 1 & 0\\
1 & 0 & 0 & 0 & 1\\\hline
N & 1 & 0 & 0 & 0\\
N & 0 & pred_{N,N} & 0 & 0\\
1 & 0 & pred_{N,1} & 0 & 0\\\hline
\end{array}
\]
where $0$ denotes the empty span and $1$ denotes the identity span. The span
$pred_{N,1}$ is the partial function which returns error on zero, and the span
$pred_{N,N}$ is the partial function returning $n-1$ for $n>0$.\bigskip
\end{ex}

We describe next a derived operation which is a minor modification of the
parallel composition, in order to simplify the matrix version of the parallel
composition. The mathematical fact behind the derived operation is this:\ in a
symmetric monoidal category with direct sums, in which the tensor product
distributes over the direct sums, if two arrows are represented as matrices,
then via distributivity isomorphisms the matrix of the tensor product of two
arrows is a tensor product of the matrices of the arrows. The precise
distributivity isomorphism needs to be specified since there are many
possible, resulting in different ordering of the rows and columns of the
tensor product matrix.

\begin{defn} \emph{Distributed parallel.}

Given systems $G_{Y_{1}+\cdots+Y_{l};A,B}^{X_{1}+\cdots+X_{k}}$,
$H_{W_{1}+\cdots+W_{l^{\prime}};B,C}^{Z_{1}+\cdots+Z_{k^{\prime}}}$ the
parallel composite $G||H$ has left interface $A$, right interface $C$, top
interface $(X_{1}+\cdots+X_{k})\times(Z_{1}+\cdots+Z_{k^{\prime}})$ and bottom
interface $(Y_{1}+\cdots+Y_{l})\times(W_{1}+\cdots+W_{l^{\prime}})$. Composing
on the top and bottom interfaces with distributivity isomorphisms we obtain a
system with left interface $A$, right interface $C$, top interface
$X_{1}\times Z_{1}+X_{2}\times Z_{1}+\cdots+X_{k}\times Z_{k^{\prime}}$ and
bottom interface $Y_{1}\times W_{1}+Y_{2}\times W_{1}+\cdots+Y_{l}\times
W_{l^{\prime}}$ The set of states of $G||H$ may similarly be distributed to
have the form $U_{1}\times V_{1}+\cdots+U_{m}\times V_{m^{\prime}}$. We will,
with an abuse of notation, denote this resulting system also as $G||H$. 
\end{defn}

\begin{defn} \emph{Distributed product.}

Given systems $G_{Y_{1}+\cdots+Y_{l};A,B}^{X_{1}+\cdots+X_{k}}$,
$H_{W_{1}+\cdots+W_{l^{\prime}};C,D}^{Z_{1}+\cdots+Z_{k^{\prime}}}$ the
product  $G\times H$ has left interface $A\times C$, right interface $B\times D$, top
interface $(X_{1}+\cdots+X_{k})\times(Z_{1}+\cdots+Z_{k^{\prime}})$ and bottom
interface $(Y_{1}+\cdots+Y_{l})\times(W_{1}+\cdots+W_{l^{\prime}})$. Composing
on the top and bottom interfaces with distributivity isomorphisms we obtain a
system with left interface $A\times C$, right interface $B\times D$, top interface
$X_{1}\times Z_{1}+X_{2}\times Z_{1}+\cdots+X_{k}\times Z_{k^{\prime}}$ and
bottom interface $Y_{1}\times W_{1}+Y_{2}\times W_{1}+\cdots+Y_{l}\times
W_{l^{\prime}}$ The set of states of $G||H$ may similarly be distributed to
have the form $U_{1}\times V_{1}+\cdots+U_{m}\times V_{m^{\prime}}$. We will,
with an abuse of notation, denote this resulting system also as $G\times H$. 
\end{defn}
\bigskip

The last assumption we make has the consequence that the pushout in sequential
composition is done a the level of control, not of data, and is therefore implementable.\bigskip

\noindent\textbf{Assumption 4. } \textit{We assume that in the matrix of the
system }$G_{B;X,Y}^{A}$\textit{ that the entries involving the sequential
interfaces are either the identity span }$1$\textit{ or the empty span }$0$.\bigskip

\subsubsection{Automaton representation}

Of course the matrix for $G_{Y,;A,B}^{X}$ has a geometric representation as
\emph{a labelled automaton}, with top sequential interfaces $X_{1}%
,X_{2},\cdots,X_{k}$, bottom sequential interfaces $Y_{1},Y_{2},\cdots,Y_{l}$,
and vertices which are labelled by the sets $U_{i}$ and for each $a\in A,b\in
B$ edges from $U_{i}$ to $U_{j}$ labelled $G_{a,b,U_{i},U_{j}}$. As usual we
will omit edges labelled with empty spans. This representation has advantages
both technical and conceptual, but is less easy to typeset. We give one
example, namely the automaton representation of the predecessor system
described above, which however has trivial parallel interfaces.%
We will see further examples in section 5.
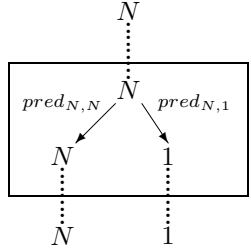
\[
\centerline{\tt\setlength{\unitlength}{2.5pt}
\begin{picture}(130,45)
\thinlines
\drawframebox{56.0}{26.0}{36.0}{20.0}{}
\drawcenteredtext{46.0}{10.0}{$N$}
\drawcenteredtext{56.0}{44.0}{$N$}
\drawcenteredtext{62.0}{10.0}{$1$}
\drawcenteredtext{56.0}{32.0}{$N$}
\drawcenteredtext{46.0}{22.0}{$N$}
\drawcenteredtext{62.0}{22.0}{$1$}
\drawvector{54.0}{30.0}{6.0}{-1}{-1}
\drawvector{58.0}{30.0}{4.0}{2}{-3}
\drawdotline{46.0}{20.0}{46.0}{12.0}
\drawdotline{62.0}{20.0}{62.0}{12.0}
\drawdotline{56.0}{42.0}{56.0}{34.0}
\drawcenteredtext{46.0}{30.0}{$\scriptstyle{pred_{N,N}}$}
\drawcenteredtext{66.0}{30.0}{$\scriptstyle{pred_{N,1}}$}
\end{picture}
}%
\]

\section{The programming language Cospan-Span}

The idea of this section is to restate the notion of system we have developed,
and describe the operations on systems. The reader should compare the notions
described here with those described in \cite{KSW00biwim} where finite state
systems were considered. We describe the programming language at the same time
as its semantics. The programs are the expressions in the operations and
constants; an execution of a program is a path in the graph described by the expression.

\subsection{Systems}

\begin{defn}
A system $G$ consists of (i) two finite sets $A,B$ called the left and right
parallel interfaces on $G;$(ii) two families of possibly infinite sets
$X=X_{1},X_{2},\cdots,X_{k\text{ }}$ and $Y=Y_{1},Y_{2},\cdots,Y_{l}$ called
the top and bottom sequential interfaces; (iii) a family of possibly infinite
sets $U=U_{1},U_{2},\cdots,U_{m}$ which together constitute the internal state
space of $G$; (iv) two functions $\varphi:\{1,2,\cdots,k\}\rightarrow
\{1,2,\cdots m\}$ and $\psi:\{1,2,\cdots,l\}\rightarrow\{1,2,\cdots m\}$
called the inclusions of the sequential interfaces, with the properties that
$X_{i}=U_{\varphi(i)}$ and $Y_{i}=U_{\psi(i)}$; (v) a family of spans of sets
$G_{a,b,i,j}:U_{i}\rightarrow U_{j}$ ($a\in A,b\in B,i\in\{1,2,\cdots
,m\},j\in\{1,2,\cdots,m\}$ which together constitute a family of graphs
$G_{a,b}$ ($a\in A,b\in B$) each with vertex set $U_{1}+U_{2}+\cdots+U_{m}.$ The graph $G_{a,b}$
is the graph of transitions of the system when the ``signals $a,b$ occur on the
parallel interfaces''. We denote the system as $G_{Y;A,B}^{X}(U)$.
\end{defn}

It is easy to see that this is the essential concrete content of the notion of
system developed in the previous section.

\subsection{Operations on systems, and constants}

In the following we denote families by giving a typical element.

\begin{defn}
The \emph{(distributed) product} of two systems $G_{Y;A,B}^{X}(U)$, $H_{W;C,D}^{Z}(V)$,
denoted $G\times H$, has left and right interfaces $A\times C,B\times D$, top
interface $\{X_{i}\times Z_{j}\},$ bottom interface $\{Y_{i}\times W_{j}\}$,
internal state space $\{U_{i}\times V_{j}\}$, inclusions of sequential
interfaces $\varphi_{G\times H}=\varphi_{G}\times\varphi_{H}$ and
$\psi_{G\times H}=\psi_{G}\times\psi_{H}$, and finally the spans
\[
(G\times H)_{(a,c),(b,d),(i_{1},j_{1}),(i_{2},j_{2})}=G_{a,,b,i_{1},i_{2}%
}\times H_{c,d,j_{1},j_{2}}.
\]
\end{defn}
Ignoring the sequential interfaces, the matrix of the distributed product is just the tensor
product of the matrices of the components.

\begin{defn}
The \emph{parallel composition} of two systems $G_{Y;A,B}^{X}(U)$,
$H_{W;B,C}^{Z}(V)$, denoted $G||H$, has left and right interfaces $A,C$, top
interface $\{X_{i}\times Z_{j}\},$ bottom interface $\{Y_{i}\times W_{j}\}$,
internal state space $\{U_{i}\times V_{j}\}$, inclusions of sequential
interfaces $\varphi_{G\times H}=\varphi_{G}\times\varphi_{H}$ and
$\psi_{G\times H}=\psi_{G}\times\psi_{H}$, and finally the spans
\[
(G||H)_{a,,c,(i_{1},j_{1}),(i_{2},j_{2})}=\sum_{b}(G_{a,b,i_{1},i_{2}}\times
H_{b,c,j_{1},j_{2}}).
\]
\end{defn}

\begin{defn}
The \emph{sequential composite} of two systems $G_{Y;A,B}^{X}(U)$, $H_{Z;C,D}^{Y}(V)$,
denoted $G\circ H$, has left and right interfaces $A+C,B+D$, top interface
$\{X_{i}\},$ bottom interface $\{Z_{i}\}$, internal state space ($\{U_{i}%
\}+\{V_{j}\}){\LARGE /}(U_{\psi_{G}(i)}\sim V_{\varphi_{H}(i)})$, inclusions
of sequential interfaces $\varphi_{G}$ and $\psi_{H}$, and finally the spans
\[
(G\circ H)_{p,q,[W_i],[W_j]}=\sum_{U\in [W_i],U'\in [W_j]}G_{p,q,U,U'}+\sum_{V\in [W_i],V'\in [W_j]}H_{p,q,V,V'}%
\]
where $p\in A+C, q\in B+D$, $W,W'\in \{U_i\}+\{V_j\}$, $[W ]$ denotes the equivalence class of $W$.
\end{defn}

\begin{defn}
The \emph{local sequential} of two systems $G_{Y;A,B}^{X}(U)$, $H_{Z;A,B}%
^{Y}(V)$, denoted $G\bullet H$, has left and right interfaces $A,B$, top
interface $\{X_{i}\},$ bottom interface $\{Z_{i}\}$, internal state space
$\{U_{i}\}+\{V_{j}\}{\LARGE /}(U_{\psi_{G}(i)}\sim V_{\varphi_{H}(i)})$, inclusions of sequential interfaces $\varphi
_{G}$ and $\psi_{H}$, and finally the spans
\[
(G\bullet H)_{p,q,[W_i],[W_j]}=\sum_{U\in [W_i],U'\in [W_j]}G_{p,q,U,U'}+\sum_{V\in [W_i],V'\in [W_j]}H_{p,q,V,V'}%
\]
where $p\in A, q\in B$, $W,W'\in \{U_i\}+\{V_j\}$, and $[W ]$ denotes the equivalence class of $W$.
\end{defn}

\begin{defn}
The\emph{ local sum} of two systems $G_{Y;A,B}^{X}(U)$, $H_{W;A,B}^{Z}(V)$,
denoted $G+ H$, has left and right interfaces $A,B$, top interface
$\{X_{i}\}+\{Z_j\},$ bottom interface $\{Y_{i}\}+\{W_{j}\}$, internal state space ($\{U_{i}%
\}+\{V_{j}\})$, inclusions of sequential interfaces $\varphi_{G}+\varphi_{H}$ and
$\psi_{G}+\psi_{H}$, and finally the spans
\[
(G+ H)_{a,c,[U_i],[U_j]}=G_{p,q,U_i,U_j}%
\]
and
\[
(G+ H)_{b,d,[V_i],[V_j]}=G_{b,d,V_,V_j},%
\]
and all remaining spans are empty.

\end{defn}

\subsection{Programs}

In our view programming languages should be presented by first describing an
algebra of systems. Then programs are elements of the free algebra of the same
type, generated by some basic systems. The meaning of the program is then
the evaluation in the concrete algebra. The programs of the Cospan-Span
language are expressions in the operations and constants of the algebra
described above, and the following basic systems: $pred_{N+1}^{N}$,
$succ_{N}^{N+1}$ (defined similarly to $pred_{N+1}^{N}$). The evaluation of a
program is a system; a behaviour is a path in the central graph of the system.

\section{Concluding remarks}

We intend in later papers to fill out details of matters sketched here, but in
fact, if one examines the previous investigations in this project it will be
clear that many matters discussed at the level of finite state control may now
be lifted to include also data.

\subsection{Turing completeness}

It is not difficult to relate the Elgot automata introduced in
\cite{W89mybook}, \cite{KSW97aprocesses},\cite{KSW00aperugia} to the algebra
of cospans of graphs. It was shown in \cite{SVW96paper with vigna} that Elgot
automata based on the elementary operations of predecessor and successor for
natural number are Turing complete, and hence also the algebra of this paper.
We give an example which illustrates sequential programming in Cospan-Span.
All the systems in the following have trivial parallel interface. In the
following we use the following constants definable from distributive category
operations, considered as systems with trivial parallel interface in which the
centre graph has no transitions (in which case a system reduces to a span of sets):
 $\eta_{X}=0\rightarrow X\overset{\nabla
}{\leftarrow}X+X,$ $\varepsilon_{X}=X+X\overset{\nabla}{\rightarrow
}X\leftarrow0$, $\nabla_{X}=X+X\overset{\nabla}{\rightarrow}X\leftarrow X,$
$1_{X}=X\overset{1}{\rightarrow}X\overset{1}{\leftarrow}X$.

\begin{ex}
The following is a program which, commencing in a state of the top sequential
interface, computes addition of two natural numbers, terminating in the lower
interface:%
\[
(\eta_{N^{2}}+1_{N^{2}})\bullet(1_{N^{2}}+\nabla)\bullet(1+pred\times
1_{N})\bullet(1_{N^{2}}+1_{N}\times succ+1_N)\bullet(\varepsilon_{N^{2}%
}+1_{N}).
\]
The system described by the program is:
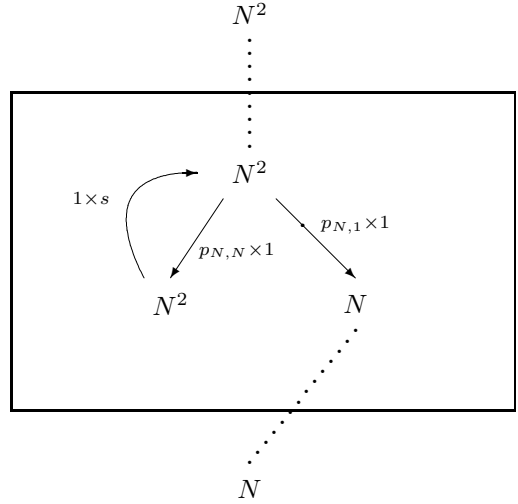
\[
\centerline{\tt\setlength{\unitlength}{5pt}
\begin{picture}(46,48)
\thinlines
\drawframebox{23.0}{24.0}{38.0}{24.0}{}
\drawcenteredtext{22.0}{42.0}{$N^2$}
\drawcenteredtext{22.0}{30.0}{$N^2$}
\drawcenteredtext{30.0}{20.0}{$N$}
\drawcenteredtext{16.0}{20.0}{$N^2$}
\drawcenteredtext{22.0}{6.0}{$N$}
\drawdotline{22.0}{40.0}{22.0}{32.0}
\drawdotline{30.0}{18.0}{22.0}{8.0}
\drawvector{24.0}{28.0}{6.0}{1}{-1}
\drawdotline{26.0}{26.0}{26.0}{26.0}
\drawvector{20.0}{28.0}{4.0}{-2}{-3}
\path
(14.0,22.0)(14.0,22.0)(13.92,22.15)(13.84,22.31)(13.77,22.47)(13.69,22.61)(13.62,22.77)(13.55,22.93)(13.48,23.08)(13.43,23.22)
\path
(13.43,23.22)(13.37,23.36)(13.3,23.52)(13.26,23.65)(13.2,23.79)(13.15,23.93)(13.11,24.08)(13.05,24.22)(13.02,24.34)(12.97,24.47)
\path
(12.97,24.47)(12.94,24.61)(12.9,24.75)(12.87,24.86)(12.84,25.0)(12.81,25.13)(12.79,25.25)(12.77,25.36)(12.75,25.5)(12.72,25.61)
\path
(12.72,25.61)(12.7,25.72)(12.7,25.84)(12.68,25.95)(12.68,26.08)(12.67,26.18)(12.65,26.29)(12.65,26.4)(12.65,26.5)(12.67,26.61)
\path
(12.67,26.61)(12.67,26.72)(12.68,26.81)(12.69,26.91)(12.7,27.02)(12.7,27.11)(12.72,27.2)(12.75,27.29)(12.77,27.4)(12.79,27.49)
\path
(12.79,27.49)(12.81,27.56)(12.85,27.65)(12.88,27.75)(12.92,27.83)(12.95,27.9)(12.98,28.0)(13.04,28.06)(13.07,28.15)(13.12,28.22)
\path
(13.12,28.22)(13.17,28.29)(13.21,28.36)(13.28,28.45)(13.32,28.52)(13.38,28.58)(13.45,28.65)(13.51,28.7)(13.57,28.77)(13.64,28.84)
\path
(13.64,28.84)(13.71,28.9)(13.79,28.95)(13.86,29.0)(13.94,29.06)(14.02,29.11)(14.1,29.18)(14.19,29.22)(14.27,29.27)(14.36,29.31)
\path
(14.36,29.31)(14.45,29.36)(14.54,29.4)(14.64,29.45)(14.73,29.5)(14.85,29.52)(14.95,29.56)(15.05,29.61)(15.15,29.63)(15.27,29.68)
\path
(15.27,29.68)(15.38,29.7)(15.5,29.74)(15.62,29.75)(15.73,29.79)(15.86,29.81)(15.98,29.84)(16.12,29.86)(16.25,29.88)(16.37,29.9)
\path
(16.37,29.9)(16.51,29.9)(16.64,29.93)(16.79,29.93)(16.93,29.95)(17.07,29.97)(17.21,29.97)(17.37,29.97)(17.53,29.99)(17.68,29.99)
\path(17.68,29.99)(17.84,29.99)(17.98,30.0)(18.0,30.0)
\drawvector{16.93}{29.95}{1.06}{1}{0}
\drawcenteredtext{21.0}{24.0}{$\scriptstyle{p_{N,N}\times1}$}
\drawcenteredtext{30.0}{26.0}{$\scriptstyle{p_{N,1}\times1}$}
\drawcenteredtext{10.0}{28.0}{$\scriptstyle{1\times s}$}
\end{picture}
}%
\]
where $p_{N,1},p_{N,N}$ are the partial functions arising from
$predecessor:N\rightarrow1+N$, and $s$ is the successor function.
\end{ex}

\subsection{Classical problems of concurrency}

We have described elsewhere (\cite{KSW00biwim},\cite{KSW00aperugia}, \cite{KSW97bspan graph}) how in $Span(Graph)$ classical problems of
concurrency may be modelled, at the level of finite state abstraction, which
is the appropriate level for controlling many properties. The current work
shows how these descriptions may be extended to include also operations on the
data types. 

We give a simple example of a parallel composite of two systems
$P$ and $Q.$ $P$ has trivial left interface, and right interface
$\{\epsilon,a\}$ whereas $Q$ has trivial right interface and left interface
$\{\epsilon,a\}$. The combined system may be represented by the diagram
(analogous to those \cite{KSW97bspan graph}), in which the first part of a 
label is the span of sets, and the
second part is the label on the parallel interface. The left system is $P$ and the right $Q$. 

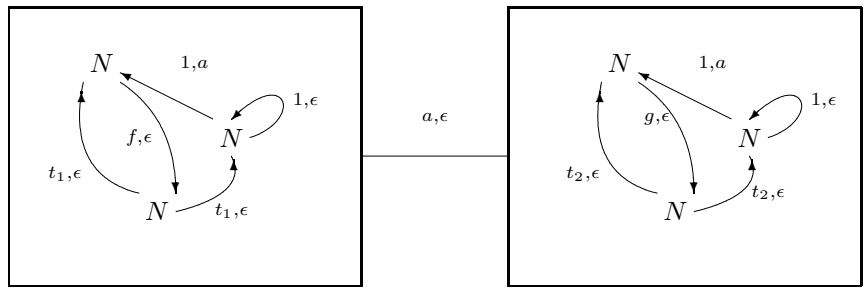
$$
\centerline{\tt\setlength{\unitlength}{3.5pt}
\begin{picture}(100,38)
\thinlines
\drawframebox{23.0}{19.0}{38.0}{30.0}{}
\drawframebox{77.0}{19.0}{38.0}{30.0}{}
\drawpath{42.0}{18.0}{58.0}{18.0}
\drawcenteredtext{50.0}{22.0}{$\scriptstyle{a,\epsilon}$}
\drawcenteredtext{14.0}{28.0}{$N$}
\drawcenteredtext{20.0}{12.0}{$N$}
\drawcenteredtext{28.0}{20.0}{$N$}
\path(16.0,26.0)(16.0,26.0)(16.11,25.9)(16.22,25.83)(16.34,25.75)(16.47,25.66)(16.58,25.59)(16.68,25.5)(16.81,25.41)(16.91,25.33)
\path(16.91,25.33)(17.02,25.24)(17.13,25.15)(17.24,25.06)(17.34,24.97)(17.45,24.88)(17.56,24.79)(17.65,24.7)(17.75,24.61)(17.86,24.52)
\path(17.86,24.52)(17.95,24.43)(18.06,24.33)(18.15,24.24)(18.25,24.13)(18.34,24.04)(18.43,23.93)(18.52,23.84)(18.61,23.75)(18.7,23.63)
\path(18.7,23.63)(18.79,23.54)(18.88,23.43)(18.97,23.34)(19.04,23.22)(19.13,23.13)(19.22,23.02)(19.29,22.91)(19.38,22.81)(19.45,22.7)
\path(19.45,22.7)(19.54,22.59)(19.61,22.49)(19.68,22.38)(19.75,22.27)(19.83,22.15)(19.9,22.04)(19.97,21.93)(20.04,21.81)(20.11,21.7)
\path(20.11,21.7)(20.18,21.59)(20.25,21.47)(20.31,21.34)(20.36,21.22)(20.43,21.11)(20.5,21.0)(20.54,20.86)(20.61,20.75)(20.66,20.63)
\path(20.66,20.63)(20.72,20.5)(20.77,20.38)(20.83,20.25)(20.88,20.13)(20.93,20.0)(20.99,19.88)(21.02,19.75)(21.08,19.63)(21.13,19.5)
\path(21.13,19.5)(21.16,19.36)(21.22,19.24)(21.25,19.11)(21.29,18.97)(21.34,18.84)(21.38,18.7)(21.41,18.56)(21.45,18.43)(21.49,18.29)
\path(21.49,18.29)(21.52,18.15)(21.56,18.02)(21.59,17.88)(21.61,17.75)(21.65,17.59)(21.68,17.45)(21.7,17.31)(21.72,17.18)(21.75,17.04)
\path(21.75,17.04)(21.77,16.88)(21.79,16.75)(21.81,16.59)(21.84,16.45)(21.86,16.31)(21.88,16.15)(21.88,16.0)(21.9,15.86)(21.91,15.7)
\path(21.91,15.7)(21.93,15.56)(21.95,15.4)(21.95,15.25)(21.97,15.09)(21.97,14.93)(21.97,14.79)(21.99,14.63)(21.99,14.47)(21.99,14.31)
\path(21.99,14.31)(21.99,14.15)(22.0,14.0)(22.0,14.0)
\drawvector{21.97}{15.09}{1.09}{0}{-1}
\path(18.0,14.0)(18.0,14.0)(17.84,14.04)(17.68,14.08)(17.52,14.11)(17.36,14.16)(17.22,14.2)(17.06,14.25)(16.91,14.31)(16.77,14.36)
\path(16.77,14.36)(16.63,14.41)(16.5,14.47)(16.36,14.52)(16.22,14.59)(16.08,14.65)(15.95,14.7)(15.81,14.77)(15.68,14.84)(15.56,14.9)
\path(15.56,14.9)(15.43,14.97)(15.31,15.04)(15.2,15.11)(15.08,15.18)(14.95,15.25)(14.84,15.34)(14.72,15.41)(14.61,15.5)(14.5,15.58)
\path(14.5,15.58)(14.4,15.65)(14.29,15.74)(14.2,15.83)(14.09,15.9)(14.0,16.0)(13.9,16.09)(13.79,16.18)(13.7,16.27)(13.61,16.38)
\path(13.61,16.38)(13.52,16.47)(13.43,16.56)(13.36,16.66)(13.27,16.77)(13.2,16.86)(13.11,16.97)(13.04,17.09)(12.95,17.18)(12.88,17.29)
\path(12.88,17.29)(12.81,17.4)(12.75,17.52)(12.68,17.63)(12.61,17.75)(12.56,17.88)(12.5,17.99)(12.43,18.11)(12.38,18.24)(12.31,18.36)
\path(12.31,18.36)(12.27,18.49)(12.22,18.61)(12.16,18.74)(12.11,18.86)(12.08,19.0)(12.04,19.13)(12.0,19.27)(11.95,19.4)(11.91,19.54)
\path(11.91,19.54)(11.88,19.68)(11.84,19.83)(11.81,19.97)(11.79,20.11)(11.75,20.27)(11.74,20.4)(11.72,20.56)(11.7,20.7)(11.68,20.86)
\path(11.68,20.86)(11.65,21.02)(11.63,21.18)(11.63,21.34)(11.61,21.49)(11.61,21.65)(11.59,21.81)(11.59,21.97)(11.59,22.15)(11.59,22.31)
\path(11.59,22.31)(11.59,22.47)(11.59,22.65)(11.59,22.83)(11.61,23.0)(11.61,23.16)(11.63,23.34)(11.63,23.52)(11.65,23.7)(11.68,23.88)
\path(11.68,23.88)(11.7,24.06)(11.72,24.25)(11.74,24.45)(11.75,24.63)(11.79,24.81)(11.81,25.0)(11.84,25.2)(11.88,25.4)(11.91,25.59)
\path(11.91,25.59)(11.95,25.79)(12.0,25.99)(12.0,26.0)
\drawvector{11.79}{24.81}{0.2}{0}{1}
\path(22.0,12.0)(22.0,12.0)(22.15,12.04)(22.31,12.08)(22.47,12.11)(22.61,12.15)(22.77,12.2)(22.91,12.24)(23.06,12.27)(23.2,12.33)
\path(23.2,12.33)(23.34,12.36)(23.5,12.41)(23.63,12.45)(23.77,12.5)(23.9,12.54)(24.04,12.59)(24.16,12.63)(24.29,12.68)(24.43,12.72)
\path(24.43,12.72)(24.54,12.77)(24.66,12.83)(24.79,12.88)(24.9,12.91)(25.02,12.97)(25.15,13.02)(25.25,13.06)(25.36,13.11)(25.47,13.16)
\path(25.47,13.16)(25.59,13.22)(25.68,13.27)(25.79,13.31)(25.9,13.36)(25.99,13.43)(26.09,13.47)(26.18,13.52)(26.27,13.59)(26.36,13.63)
\path(26.36,13.63)(26.45,13.68)(26.54,13.75)(26.63,13.79)(26.7,13.86)(26.79,13.9)(26.86,13.97)(26.95,14.02)(27.02,14.08)(27.09,14.13)
\path(27.09,14.13)(27.16,14.2)(27.24,14.25)(27.31,14.31)(27.36,14.38)(27.43,14.43)(27.5,14.5)(27.54,14.56)(27.61,14.61)(27.66,14.68)
\path(27.66,14.68)(27.72,14.74)(27.77,14.79)(27.81,14.86)(27.86,14.91)(27.9,14.99)(27.95,15.04)(28.0,15.11)(28.02,15.18)(28.06,15.24)
\path(28.06,15.24)(28.11,15.31)(28.13,15.36)(28.16,15.43)(28.2,15.5)(28.22,15.56)(28.25,15.63)(28.27,15.7)(28.29,15.77)(28.31,15.84)
\path(28.31,15.84)(28.33,15.9)(28.34,15.97)(28.36,16.04)(28.36,16.11)(28.38,16.18)(28.38,16.25)(28.38,16.33)(28.38,16.4)(28.4,16.47)
\path(28.4,16.47)(28.38,16.54)(28.38,16.61)(28.38,16.68)(28.38,16.77)(28.36,16.84)(28.36,16.9)(28.34,16.99)(28.33,17.06)(28.31,17.13)
\path(28.31,17.13)(28.29,17.2)(28.27,17.29)(28.25,17.36)(28.22,17.43)(28.2,17.52)(28.16,17.59)(28.13,17.68)(28.11,17.75)(28.06,17.84)
\path(28.06,17.84)(28.02,17.91)(28.0,17.99)(28.0,18.0)
\drawvector{28.27}{17.29}{0.27}{0}{1}
\path(30.0,20.0)(30.0,20.0)(30.16,20.06)(30.34,20.11)(30.5,20.18)(30.68,20.25)(30.84,20.31)(30.99,20.38)(31.13,20.47)(31.27,20.54)
\path(31.27,20.54)(31.41,20.61)(31.56,20.7)(31.68,20.77)(31.81,20.86)(31.93,20.93)(32.05,21.02)(32.16,21.11)(32.27,21.2)(32.37,21.27)
\path(32.37,21.27)(32.48,21.36)(32.58,21.45)(32.66,21.54)(32.75,21.63)(32.83,21.72)(32.91,21.81)(32.98,21.9)(33.05,22.0)(33.12,22.08)
\path(33.12,22.08)(33.18,22.16)(33.23,22.25)(33.29,22.34)(33.33,22.43)(33.38,22.52)(33.43,22.61)(33.47,22.7)(33.5,22.79)(33.52,22.88)
\path(33.52,22.88)(33.55,22.95)(33.56,23.04)(33.58,23.13)(33.59,23.2)(33.61,23.29)(33.62,23.36)(33.62,23.45)(33.62,23.52)(33.61,23.59)
\path(33.61,23.59)(33.59,23.66)(33.58,23.74)(33.56,23.79)(33.55,23.86)(33.51,23.93)(33.5,24.0)(33.45,24.04)(33.43,24.11)(33.38,24.15)
\path(33.38,24.15)(33.34,24.2)(33.3,24.25)(33.25,24.31)(33.19,24.34)(33.13,24.38)(33.08,24.43)(33.01,24.45)(32.94,24.49)(32.87,24.5)
\path(32.87,24.5)(32.8,24.54)(32.73,24.56)(32.65,24.56)(32.56,24.58)(32.48,24.59)(32.4,24.59)(32.3,24.59)(32.2,24.59)(32.11,24.58)
\path(32.11,24.58)(32.01,24.56)(31.9,24.54)(31.79,24.52)(31.68,24.5)(31.56,24.45)(31.45,24.43)(31.33,24.38)(31.2,24.34)(31.08,24.27)
\path(31.08,24.27)(30.95,24.22)(30.81,24.15)(30.68,24.09)(30.54,24.02)(30.4,23.93)(30.27,23.84)(30.11,23.75)(29.97,23.65)(29.81,23.56)
\path(29.81,23.56)(29.66,23.45)(29.5,23.34)(29.36,23.2)(29.18,23.08)(29.02,22.95)(28.86,22.81)(28.7,22.65)(28.52,22.5)(28.34,22.34)
\path(28.34,22.34)(28.16,22.16)(28.0,22.0)(28.0,22.0)
\drawvector{29.02}{22.95}{1.02}{-1}{-1}
\drawvector{26.0}{22.0}{10.0}{-2}{1}
\drawcenteredtext{70.0}{28.0}{$N$}
\drawcenteredtext{76.0}{12.0}{$N$}
\drawcenteredtext{84.0}{20.0}{$N$}
\path(72.0,26.0)(72.0,26.0)(72.11,25.9)(72.22,25.83)(72.34,25.75)(72.47,25.66)(72.58,25.59)(72.69,25.5)(72.8,25.41)(72.91,25.33)
\path(72.91,25.33)(73.02,25.24)(73.13,25.15)(73.24,25.06)(73.34,24.97)(73.44,24.88)(73.55,24.79)(73.66,24.7)(73.75,24.61)(73.86,24.52)
\path(73.86,24.52)(73.95,24.43)(74.05,24.33)(74.15,24.24)(74.25,24.13)(74.33,24.04)(74.44,23.93)(74.52,23.84)(74.61,23.75)(74.7,23.63)
\path(74.7,23.63)(74.8,23.54)(74.88,23.43)(74.97,23.34)(75.05,23.22)(75.13,23.13)(75.22,23.02)(75.3,22.91)(75.38,22.81)(75.45,22.7)
\path(75.45,22.7)(75.54,22.59)(75.61,22.49)(75.69,22.38)(75.75,22.27)(75.83,22.15)(75.91,22.04)(75.97,21.93)(76.05,21.81)(76.11,21.7)
\path(76.11,21.7)(76.18,21.59)(76.25,21.47)(76.3,21.34)(76.36,21.22)(76.43,21.11)(76.5,21.0)(76.55,20.86)(76.61,20.75)(76.66,20.63)
\path(76.66,20.63)(76.72,20.5)(76.77,20.38)(76.83,20.25)(76.88,20.13)(76.94,20.0)(76.99,19.88)(77.02,19.75)(77.08,19.63)(77.13,19.5)
\path(77.13,19.5)(77.16,19.36)(77.22,19.24)(77.25,19.11)(77.3,18.97)(77.33,18.84)(77.38,18.7)(77.41,18.56)(77.44,18.43)(77.49,18.29)
\path(77.49,18.29)(77.52,18.15)(77.55,18.02)(77.58,17.88)(77.61,17.75)(77.65,17.59)(77.68,17.45)(77.69,17.31)(77.72,17.18)(77.75,17.04)
\path(77.75,17.04)(77.77,16.88)(77.8,16.75)(77.81,16.59)(77.83,16.45)(77.86,16.31)(77.88,16.15)(77.88,16.0)(77.91,15.86)(77.91,15.7)
\path(77.91,15.7)(77.93,15.56)(77.94,15.4)(77.95,15.25)(77.97,15.09)(77.97,14.93)(77.97,14.79)(77.99,14.63)(77.99,14.47)(77.99,14.31)
\path(77.99,14.31)(77.99,14.15)(78.0,14.0)(78.0,14.0)
\drawvector{77.97}{15.09}{1.09}{0}{-1}
\path(74.0,14.0)(74.0,14.0)(73.83,14.04)(73.68,14.08)(73.52,14.11)(73.36,14.16)(73.22,14.2)(73.06,14.25)(72.91,14.31)(72.77,14.36)
\path(72.77,14.36)(72.63,14.41)(72.5,14.47)(72.36,14.52)(72.22,14.59)(72.08,14.65)(71.94,14.7)(71.81,14.77)(71.69,14.84)(71.55,14.9)
\path(71.55,14.9)(71.44,14.97)(71.31,15.04)(71.19,15.11)(71.08,15.18)(70.95,15.25)(70.83,15.34)(70.72,15.41)(70.61,15.5)(70.5,15.58)
\path(70.5,15.58)(70.4,15.65)(70.3,15.74)(70.19,15.83)(70.08,15.9)(70.0,16.0)(69.9,16.09)(69.8,16.18)(69.7,16.27)(69.61,16.38)
\path(69.61,16.38)(69.52,16.47)(69.44,16.56)(69.36,16.66)(69.27,16.77)(69.19,16.86)(69.11,16.97)(69.04,17.09)(68.95,17.18)(68.88,17.29)
\path(68.88,17.29)(68.81,17.4)(68.75,17.52)(68.68,17.63)(68.61,17.75)(68.55,17.88)(68.5,17.99)(68.44,18.11)(68.38,18.24)(68.31,18.36)
\path(68.31,18.36)(68.27,18.49)(68.22,18.61)(68.16,18.74)(68.11,18.86)(68.08,19.0)(68.04,19.13)(68.0,19.27)(67.95,19.4)(67.91,19.54)
\path(67.91,19.54)(67.88,19.68)(67.84,19.83)(67.81,19.97)(67.79,20.11)(67.75,20.27)(67.74,20.4)(67.72,20.56)(67.69,20.7)(67.68,20.86)
\path(67.68,20.86)(67.66,21.02)(67.63,21.18)(67.63,21.34)(67.61,21.49)(67.61,21.65)(67.59,21.81)(67.59,21.97)(67.59,22.15)(67.58,22.31)
\path(67.58,22.31)(67.59,22.47)(67.59,22.65)(67.59,22.83)(67.61,23.0)(67.61,23.16)(67.63,23.34)(67.63,23.52)(67.66,23.7)(67.68,23.88)
\path(67.68,23.88)(67.69,24.06)(67.72,24.25)(67.74,24.45)(67.75,24.63)(67.79,24.81)(67.81,25.0)(67.84,25.2)(67.88,25.4)(67.91,25.59)
\path(67.91,25.59)(67.95,25.79)(68.0,25.99)(68.0,26.0)
\drawvector{67.79}{24.81}{0.2}{0}{1}
\path(78.0,12.0)(78.0,12.0)(78.15,12.04)(78.3,12.08)(78.47,12.11)(78.61,12.15)(78.77,12.2)(78.91,12.24)(79.06,12.27)(79.2,12.33)
\path(79.2,12.33)(79.34,12.36)(79.5,12.41)(79.63,12.45)(79.77,12.5)(79.91,12.54)(80.04,12.59)(80.16,12.63)(80.3,12.68)(80.43,12.72)
\path(80.43,12.72)(80.55,12.77)(80.66,12.83)(80.79,12.88)(80.91,12.91)(81.02,12.97)(81.15,13.02)(81.25,13.06)(81.36,13.11)(81.47,13.16)
\path(81.47,13.16)(81.58,13.22)(81.69,13.27)(81.79,13.31)(81.9,13.36)(81.99,13.43)(82.08,13.47)(82.19,13.52)(82.27,13.59)(82.36,13.63)
\path(82.36,13.63)(82.45,13.68)(82.55,13.75)(82.63,13.79)(82.7,13.86)(82.79,13.9)(82.86,13.97)(82.94,14.02)(83.02,14.08)(83.09,14.13)
\path(83.09,14.13)(83.16,14.2)(83.24,14.25)(83.3,14.31)(83.36,14.38)(83.43,14.43)(83.5,14.5)(83.55,14.56)(83.61,14.61)(83.66,14.68)
\path(83.66,14.68)(83.72,14.74)(83.77,14.79)(83.81,14.86)(83.86,14.91)(83.91,14.99)(83.94,15.04)(84.0,15.11)(84.02,15.18)(84.06,15.24)
\path(84.06,15.24)(84.11,15.31)(84.13,15.36)(84.16,15.43)(84.19,15.5)(84.22,15.56)(84.25,15.63)(84.27,15.7)(84.29,15.77)(84.3,15.84)
\path(84.3,15.84)(84.33,15.9)(84.34,15.97)(84.36,16.04)(84.36,16.11)(84.38,16.18)(84.38,16.25)(84.38,16.33)(84.38,16.4)(84.4,16.47)
\path(84.4,16.47)(84.38,16.54)(84.38,16.61)(84.38,16.68)(84.38,16.77)(84.36,16.84)(84.36,16.9)(84.34,16.99)(84.33,17.06)(84.3,17.13)
\path(84.3,17.13)(84.29,17.2)(84.27,17.29)(84.25,17.36)(84.22,17.43)(84.19,17.52)(84.16,17.59)(84.13,17.68)(84.11,17.75)(84.06,17.84)
\path(84.06,17.84)(84.02,17.91)(84.0,17.99)(84.0,18.0)
\drawvector{84.27}{17.29}{0.27}{0}{1}
\path(86.0,20.0)(86.0,20.0)(86.16,20.06)(86.34,20.11)(86.5,20.18)(86.68,20.25)(86.83,20.31)(86.99,20.38)(87.13,20.47)(87.27,20.54)
\path(87.27,20.54)(87.41,20.61)(87.55,20.7)(87.69,20.77)(87.81,20.86)(87.94,20.93)(88.05,21.02)(88.16,21.11)(88.27,21.2)(88.37,21.27)
\path(88.37,21.27)(88.48,21.36)(88.58,21.45)(88.66,21.54)(88.75,21.63)(88.83,21.72)(88.91,21.81)(88.98,21.9)(89.05,22.0)(89.12,22.08)
\path(89.12,22.08)(89.18,22.16)(89.23,22.25)(89.29,22.34)(89.33,22.43)(89.38,22.52)(89.43,22.61)(89.47,22.7)(89.5,22.79)(89.52,22.88)
\path(89.52,22.88)(89.55,22.95)(89.56,23.04)(89.58,23.13)(89.59,23.2)(89.61,23.29)(89.62,23.36)(89.62,23.45)(89.62,23.52)(89.61,23.59)
\path(89.61,23.59)(89.59,23.66)(89.58,23.74)(89.56,23.79)(89.55,23.86)(89.51,23.93)(89.5,24.0)(89.45,24.04)(89.43,24.11)(89.38,24.15)
\path(89.38,24.15)(89.34,24.2)(89.3,24.25)(89.25,24.31)(89.19,24.34)(89.13,24.38)(89.08,24.43)(89.01,24.45)(88.94,24.49)(88.87,24.5)
\path(88.87,24.5)(88.8,24.54)(88.73,24.56)(88.65,24.56)(88.56,24.58)(88.48,24.59)(88.4,24.59)(88.3,24.59)(88.2,24.59)(88.11,24.58)
\path(88.11,24.58)(88.01,24.56)(87.9,24.54)(87.79,24.52)(87.68,24.5)(87.56,24.45)(87.44,24.43)(87.33,24.38)(87.2,24.34)(87.08,24.27)
\path(87.08,24.27)(86.94,24.22)(86.81,24.15)(86.69,24.09)(86.55,24.02)(86.41,23.93)(86.27,23.84)(86.11,23.75)(85.97,23.65)(85.81,23.56)
\path(85.81,23.56)(85.66,23.45)(85.5,23.34)(85.36,23.2)(85.19,23.08)(85.02,22.95)(84.86,22.81)(84.69,22.65)(84.52,22.5)(84.34,22.34)
\path(84.34,22.34)(84.16,22.16)(84.0,22.0)(84.0,22.0)
\drawvector{85.02}{22.95}{1.02}{-1}{-1}
\drawvector{82.0}{22.0}{10.0}{-2}{1}
\drawcenteredtext{36.0}{24.0}{$\scriptstyle{1,\epsilon}$}
\drawcenteredtext{24.0}{28.0}{$\scriptstyle{1,a}$}
\drawcenteredtext{10.0}{16.0}{$\scriptstyle{t_1,\epsilon}$}
\drawcenteredtext{28.0}{12.0}{$\scriptstyle{t_1,\epsilon}$}
\drawcenteredtext{18.0}{20.0}{$\scriptstyle{f,\epsilon}$}
\drawcenteredtext{92.0}{24.0}{$\scriptstyle{1,\epsilon}$}
\drawcenteredtext{80.0}{28.0}{$\scriptstyle{1,a}$}
\drawcenteredtext{66.0}{16.0}{$\scriptstyle{t_2,\epsilon}$}
\drawcenteredtext{86.0}{14.0}{$\scriptstyle{t_2,\epsilon}$}
\drawcenteredtext{74.0}{22.0}{$\scriptstyle{g,\epsilon}$}
\end{picture}
}$$
The system $P$ repeatedly applies $f$ and then a test $t_1$ until the test
results false, and then $P$ may idle, eventually (in the Italian sense) 
synchronizing with $Q$ on the signal $a$. After this $P$ repeats the whole sequence.  
$Q$ does the same, but with a different function $g$ and a different test $t_2$, and seeks to synchronize with $P$.

Each of $P$ and $Q$ may be described by a Cospan-Span program in a similar way to the addition program above.
\subsection{Hierarchy}

There is an obvious relevance to hierarchical systems of the fact that systems
in this algebra may be constructed by repeated parallel and sequential operations, with analogies
to state charts.

\subsection{Change of geometry}

Already in \cite{KSW00biwim} we discussed the description of changing geometry
using sequential operations on parallel systems. However in that paper we
considered only the local sequential composition, whereas in this paper we
have a general sequential operation, which allows change of geometry with a
change of parallel protocol. In that article we abstracted away data.

\subsection{Relation with other work}

Theoretical considerations behind this work include \cite{RSW08calculating
colimits montanari},\cite{RSWCT04generic sep alg}%
,\cite{CLW93carbonilackwalters},\cite{CW87carboniwalters},\cite{CL01extensive completion} and \cite{KSW02feedback}.

Studying \cite{dFASW08dmarkov} and \cite{dFASW08equantum} the reader will note
similarities with the algebra here. In fact, this paper is the result of
comparing {\cite{KSW97bspan graph}} with \cite{dFASW08dmarkov} and
\cite{dFASW08equantum}.

\end{document}